\documentclass[a4paper,twoside,11pt]{article}  

\usepackage{amsmath,amsthm,amssymb}
\usepackage{latexsym}
\usepackage{amsfonts}
\usepackage{amsmath}
\usepackage{amssymb}
\usepackage{amsthm}
\usepackage[active]{srcltx}
\usepackage[ansinew]{inputenc}   
\usepackage{fancyvrb}
\usepackage{array}  
\usepackage{graphicx}
\usepackage{float}

\usepackage[a4paper=true,%
breaklinks=true,%
colorlinks=true,%
linkcolor={blue},%
citecolor={red},%
pdftitle={The question ``How many $1$'s are needed?''
revisited},%
pdfauthor={J. Arias de Reyna, J. van de Lune},%
pdfkeywords={natural numbers, complexity},%
pdfstartview=FitH,
]{hyperref}

\newtheorem{theorem}{Theorem}
\newtheorem{definition}[theorem]{Definition}

\newtheorem{conjecture}{Conjecture}

\newfont{\cmbsy}{cmbsy10}
\newfont{\cmmib}{cmmib10}

\newcommand{\Orden}{\mathop{\hbox{\cmbsy O}}}

\newcommand{\AddExc}{\mathop{\hbox{\rm AddExc}}}

\def\CR{\mathrm{CR\,}}

\newcommand{\seqnum}[1]{\href{http://oeis.org/eisA.cgi?Anum=#1}{\underline{#1}}}

\begin{document}

\title{The question ``How many $1$'s are needed\,?'' revisited}

\author{J. Arias de Reyna
\and J. van de Lune}

\pagestyle{myheadings}
\markboth{Arias de Reyna and van de Lune}{The question ``How many $1$'s are needed$\,$?'' revisited}

\date{February 11, 2009}

\maketitle

\begin{abstract}
We present a rigorous and relatively fast method for the computation
of the \emph{complexity}  of a natural number (sequence
\seqnum{A005245}), and answer some \emph{old and new} questions
related to the question in the title of this note.  We also extend
the known terms of the related sequence \seqnum{A005520}.
\end{abstract}

\setcounter{section}{-1}

\section*{Introduction.}\label{intro}

The subject of this note was (more or less indirectly) initiated  in
1953 by K. Mahler and J. Popken \cite{MP}. We begin with a brief
description of part of their work: Given a symbol $x$,  consider the
set $V_n$ of all formal \emph{sum-products} which can be constructed
by using only the symbol $x$ and precisely $n-1$ symbols from
$\{+,\times\}$ and an arbitrary number of parentheses ``$($'' and
``$)$''.

We have, for example, $V_1=\{x\}$, $V_2=\{x+x,x\times x\}$,
\begin{displaymath}
V_3=\{x+(x+x),x+(x\times x), x\times(x+x), x\times(x\times x)\}.
\end{displaymath}
More generally, for $n\ge2$,
\begin{displaymath}
V_n=\bigcup_{k=1}^{n-1} (V_k+V_{n-k})\cup \bigcup_{k=1}^{n-1}
(V_k\times V_{n-k}).
\end{displaymath}

Mahler and Popken's question was the following: If $x$ is a positive
real number, what is the largest number in $V_n\,$? We restrict
ourselves here to the case $x=1$. Then the answer is \cite{MP}
$M_n:=\max V_n=\max_{1\le k\le n}p_{n,k}$ where
\begin{displaymath}
p_{n,k}=\Bigl\lfloor \frac{n}{k}\Bigr\rfloor^{k\left(\left\lfloor
\frac{n}{k}\right\rfloor+1\right)-n}
\Bigl(\Bigl\lfloor\frac{n}{k}\Bigr\rfloor
+1\Bigr)^{n-k\left\lfloor\frac{n}{k}\right\rfloor}.
\end{displaymath}
This formula was simplified by Selfridge (see Guy
\cite[p.~189]{Guy3}) to $M_{3m-1}=2\cdot 3^{m-1}$, $M_{3m}=3^m$,
 $M_{3m+1}=4\cdot 3^{m-1}$ for all $m\ge1$.  Clearly $M_1=1$.

Our problem is more or less the converse: Write a given natural
number $n$ as a \emph{sum-product} as described above, only using
the five symbols $1$, $+$, $\times$, $($, and $)$. (However, not all these
signs need to be used.)

It is clear that this is always possible: $n=1+1+1+\cdots+1$ (using
$n$ $1$'s). Some further simple examples are
\begin{displaymath}
5=1+(1+1)\times(1+1),\quad 6=(1+1)\times(1+1+1).
\end{displaymath}
Our goal will, of course, be to minimize the number of $1$'s used.

In a sum-product representation of $n$ we will usually write $2$
 instead of  $1+1$, and $3$ instead of
 $1+1+1$. Also, we  will replace  the symbol $\times$ (times)
  by a dot $\cdot$ or simply juxtapose. For example,
the Fibonacci number $F_{25} $ can then be written with 35 $1$'s as
follows
\begin{multline*}
F_{25}=75025=(1+2^2)(1+2^2)(1+2(1+2^2)(1+2^2)(2^2\cdot 3(1+2^2)))
\end{multline*}
and $2^{27}-1$ can be written with $56$ as
\begin{displaymath}
2^{27}-1=134217727=(1+2\cdot3)(1+2^3\cdot
3^2)(1+2^9\cdot3^3(1+2\cdot3^2)).
\end{displaymath}
All these examples are \emph{minimal}  in the sense defined in the
next section.

\setcounter{section}{0}

\section{Definitions and first properties.}

\begin{definition}
The minimal number of $1$'s needed to represent $n$ as a sum-product
will be denoted by $\Vert n \Vert$ and will be called the
\emph{complexity} of $n$.
\end{definition}

It is clear that $\Vert1\Vert=1$ and $\Vert2\Vert=2$, but
$\Vert11\Vert\ne2$  (``pasting together'' two 1's is not an allowed
operation).  One may verify directly that $\Vert3\Vert=3$,
$\Vert4\Vert=4$, $\Vert5\Vert=5$, $\Vert6\Vert=5$, and by means of
our program in Section \ref{sectiondos} it may be shown that
\begin{gather*}
\Vert7\Vert=6,\quad \Vert8\Vert=6,\quad \Vert9\Vert=6,\quad
\Vert10\Vert=7,\quad \Vert11\Vert=8,\quad \Vert12\Vert=7,\\
\Vert13\Vert=8,\quad \Vert14\Vert=8,\quad \Vert15\Vert=8,\quad\Vert
16\Vert=8,\quad \Vert 17\Vert=9,\quad\Vert18\Vert=8.
\end{gather*}

Note that

(a) $\Vert n\Vert$ is not monotonic

(b) $n$ may have different minimal representations
($4=1+1+1+1=(1+1)(1+1)$).

It is clear that we always have
\begin{displaymath}
\Vert a+b\Vert\le\Vert a\Vert+\Vert b\Vert\quad\text{and}\quad \Vert
a\cdot b\Vert\le\Vert a\Vert+\Vert b\Vert
\end{displaymath}
so that, for example, $\Vert 2^n\Vert\le 2n$. Also see Section
\ref{section}.

Some useful bounds on the complexity are known
\begin{displaymath}
\frac{3}{\log 3}\log n\le \Vert n\Vert\le \frac{3}{\log2}\log
n,\qquad n\ge2.
\end{displaymath}
The first can be found in Guy \cite{Guy3} and is essentially due to
Selfridge. The second appeared in Arias de Reyna \cite{Arias} (this
inequality can easily be proved. Indeed, just think of the binary
expansion of $n$.) Since it is known that $\Vert 3^k\Vert=3k$ for
$k\ge1$, the first inequality cannot be improved. As for the second
one:
\begin{displaymath}\limsup_{n\to\infty} \Vert n\Vert/\log n\end{displaymath} is not
known, but we conjecture that it is considerably $<\frac{3}{\log2}$
($\approx 4.328$). Our most extreme observation is
$\Vert1439\Vert/\log1439\approx 3.575503$.

\section{Computing the complexity.}\label{sectiondos}

For $n\ge2$ we may compute $\Vert n\Vert$ from
\begin{equation}\label{defCompl}
\Vert n\Vert=\min\bigl\{\min_{1\le j\le n/2} \Vert j\Vert+\Vert
n-j\Vert\;,\quad \min_{d\mid n,\;2\le d\le\sqrt{n}} \Vert
d\Vert+\Vert n/d\Vert\bigr\}\,.
\end{equation}

  From this it is clear that, for large $n$, the computation of
\begin{displaymath}
\min_{1\le j\le n/2} \Vert j\Vert+\Vert n-j\Vert
\end{displaymath} is quite time consuming, if not eventually prohibitive.
Rawsthorne \cite[p.~14]{Rawsthorne} wrote \emph{This formula is very
time consuming to use for large $n$,  but we know of no other way to
calculate $\Vert n\Vert$. }

The principal goal of this note is to reduce  the number of
operations for the computation of $\Vert n\Vert$. (We can show that,
instead of $\Orden(n^2)$,  our algorithm  needs only
$\Orden(n^{1.345})$ operations for the computation of $\Vert
n\Vert$.)

According to the definition we have to compute
\begin{displaymath}
P:=\min_{1\le k\le n/2} \Vert k\Vert+\Vert n-k\Vert\quad \text{and}
\quad T:=\min_{d\mid n,\;2\le d\le\sqrt{n}} \Vert d\Vert+\Vert
n/d\Vert
\end{displaymath}
and then set $\Vert n\Vert =\min(P,T)$. It is clear that $P\le
\Vert1\Vert+\Vert n-1\Vert$ so that $P$ is the result of the loop
\begin{align*}
&P=1+\Vert n-1\Vert;\\
&\text{For\;\;}k=2\;\text{ to }\; k = n/2\;\; \text{ do } \quad
P=\min(P\;,\;\;\Vert k\Vert+\Vert n-k\Vert)\;.\;
\end{align*}
Clearly this is cumbersome for large $n$. It would be very helpful
to have a relatively small number $\textbf{kMax}$ such that $P$
would just as well be the result of the much shorter loop
\begin{align*}
&P=1+\Vert n-1\Vert;\\
&\text{For\;\;}k=2\;\text{ to }\; k = \textbf{kMax}\;\; \text{ do }
\quad P=\min(P\;,\;\;\Vert k\Vert+\Vert n-k\Vert)\;.\;
\end{align*}
Such a relatively small $\textbf{kMax}$ may be found indeed by observing that
\begin{displaymath}
\Vert m\Vert\ge \frac{3}{\log 3}\log m\qquad \text{for all $m\ge1$.}
\end{displaymath}
Indeed, we are through if $\textbf{kMax}$ satisfies
\begin{multline*}
\Vert k\Vert+\Vert n-k\Vert\ge \frac{3}{\log 3}\bigl(\log k+\log
(n-k)\bigr)\ge 1+\Vert n-1\Vert\\\text{  for $\textbf{kMax}+1\le
k\le n/2$.}
\end{multline*}
This only requires to solve a simple quadratic inequality:
\begin{displaymath}
k^2-nk+\exp(R)\ge0 \qquad \text{where} \quad R=\frac{\log
3}{3}(1+\Vert n-1\Vert).
\end{displaymath}
It is easily seen that, for $n\ge 7$, we can safely take
\begin{displaymath}
\textbf{kMax}=\Bigl\lfloor\frac12+
\frac{n}{2}\Bigl(1-\sqrt{1-4\exp(R-2\log n)}\Bigr)\Bigr\rfloor.
\end{displaymath}
It will soon become clear that for large $n$ this $\textbf{kMax}$ is
very small compared to $n/2$.  In our computations covering all
$n\le 905\,000\,000$  we observed that $\textbf{kMax}\le66$ in all
cases, with an average value  of about $11.57$.

However, we can not use this ``trick'' for the $\times$ part.
\bigskip
\goodbreak

\textbf{Mathematica program  to compute
\texttt{Compl[}$n$\texttt{]}$:=\Vert n\Vert$.}
\medskip
\begin{verbatim}
Compl[1] = 1; Compl[2] = 2; Compl[3] = 3; Compl[4] = 4; 
Compl[5] = 5; Compl[6] = 5; nDone = 6; 
(* Our computed kMax is not real for n<= 6 *) 
ComplChamp = 5; 
(* = largest value of  C[n]  found so far.*)

n = nDone; While[0 == 0, n += 1;

(*  First we deal with the  PLUS-part.  *) 
P = 1 + Compl[n - 1]; R=N[Log[3] P/3]; 
kMax =  Floor[1/2+n(1-Sqrt[ 1 - 4Exp[R - 2Log[n]]])/2]; 
For[k = 2 , k <= kMax , k++ , 
P = Min[ P , Compl[k] + Compl[n - k]]]; 
(* kMax < 2  causes no problem.  *)

(*  Now for the  TIMES-part.  *) 
S = Divisors[n];   LSplus1 = Length[S] + 1; T = P;      
(* From the PLUS-part we already
                           know  that  Compl[n] <= P  *)
For[k = 2 , k <= LSplus1/2 , k++ , 
d = S[[k]];  T = Min[ T , Compl[d] + Compl[n/d]]];

Compl[n] = T;   (*  There we are !  *)

(* We output the Champion Compl[n] and the
                         corresponding Compl[n] / Log[n] *)
If[T > ComplChamp,
    ComplChamp = T;
    Print["n = ", n, "    kMax = ", kMax, "    ComplChamp = ",
                 ComplChamp, "    ||n||/Log [n] = ",
                                        N[Compl[n]/Log[n]]]]]
\end{verbatim}

A much faster Delphi-Object-Pascal version of this program, run on a
Toshiba laptop, computes $\Vert n\Vert $ for all $n\le
905\,000\,000$ in about 2 hours and 40 minutes.
\medskip

\textbf{Note.} In the range $n\le 905\,000\,000$  it suffices to
take $\texttt{kMax}=6$. This value ($\texttt{kMax}=6$) is necessary
only for $n=353\, 942\, 783 $ and $ n = 516\, 743\, 639$. But this
is hindsight!

\section{Some records.}

\begin{definition} The number $n$ is called \emph{highly complex} if
$\Vert k\Vert<\Vert n\Vert$ for all $k<n$.
\end{definition}

P. Fabian (see~\cite{Pegg}) has computed the first 58 highly complex
numbers. With our new method we have been able to add those with
$59\le \Vert n\Vert\le 67$ (in boldface at the end of Table 1).
There are no others in the range $n\le 905\,000\,000$.  We performed
our computations on  a Toshiba laptop, 2GB RAM, 3.2 GHz, and  could
verify Fabian's results within 138 seconds.

We denote by $F_m$ the first number having complexity $m$
(i.~e.~$F_m$ is the $m$-th highly complex number). $S(m)$ denotes
the set of numbers with complexity $m$, its first element is $F_m$, and
its maximal element $M_m$.


{\scriptsize
\begin{center}

\text{ TABLE 1}

Some data related to Highly Complex numbers
\medskip

\begin{tabular}{|r|r|r|l|r|r|r|r|c|r|}
\hline
$m$ & $F_m$ & \textbf{kMax} & $\Vert F_m\Vert/\log F_m$ & $M_m$ & $\# S(m)$\\
\hline 
\phantom{(}1  &   1       &      &   & 1 & 1\cr 
2  &   2 &      & 2.8853900818  & 2 & 1\cr 
3  &  3        &      & 2.7307176799 & 3 & 1\cr 
4  &   4       &      & 2.8853900818  & 4 & 1\cr 
5  &   5       &      &  3.1066746728 &  6  & 2\cr 
6 &    7 &     &   3.0833900542 & 9 & 3\cr 
7  &  10        &  2    & 3.0400613733 & 12 & 2\cr 
8  &  11        &  2    &  3.3362591314 & 18 & 6\cr 
9 &  17         &  2    &  3.1766051148 & 27 & 6\cr 
10 & 22        &  2    &  3.2351545315 & 36 & 7\cr 
\hline
\end{tabular}
\end{center}

\begin{center}
\begin{tabular}{|r|r|r|l|r|r|r|l|l|r|}
\hline

\hline
\end{tabular}
\end{center}

\begin{center}
\begin{tabular}{|r|r|r|l|r|r|r|l|l|r|}
\hline
$m$ & $F_m$ & \textbf{kMax} & $\Vert F_m\Vert/\log F_m$ & $M_m$ & $\# S(m)$\\
\hline
11 &  23        &  3 &  3.5082188779 & 54 & 14\cr 
12 &  41        &  2    &  3.2313900968& 81 & 16\cr 
13 &  47        &  3    &  3.3764939282 & 108 & 20\cr
14 &  59        &  3    &  3.4334448653 & 162 & 34\cr 
15 &  89&  3    &  3.3417721474 & 243 & 42\cr
16 &  107       &  3    &  3.4240500919 & 324 & 56\cr
17 &  167       &  3    &  3.3216140197 & 486 & 84\cr
18 &  179       &  4    &  3.4699559034 & 729 & 108\cr
19 &  263       &  4    &  3.4098124155 & 972 & 152\cr
20 &  347       &  4    &  3.4191980703 & 1458 & 214\cr
21 &  467       &  5    &  3.4166734517 & 2187 & 295\cr
22 &  683       &  5    &  3.3708752513 & 2916 & 398\cr
23 &  719       &  6    &  3.4965771927 & 4374 & 569\cr
24 &  1223      &  5    &  3.3759727432 & 6561 & 763\cr
25 &  1438      &  7    &  3.4383125626 & 8748 & 1094\cr
26 &  1439      &  10  &  \textbf{\scriptsize3.5755032174} & 13122 & 1475\cr
27 &  2879      &  7    &  3.3897461199 & 19683 & 2058\cr
28 &  3767      &  8    &  3.4005202424 & 26244 & 2878\cr
29 &  4283      &  10  &   3.4679002280 & 39366 & 3929\cr
30 &  6299      &  9    &  3.4292979813 & 59049 & 5493\cr
31 &  10079     &  8    &  3.3629090954 & 78732 & 7669\cr
32 &  11807     &  10  &   3.4128062668 & 118098 & 10501\cr
33 &  15287     &  12  &   3.4250989750 & 177147 & 14707\cr
34 &  21599     &  12  &   3.4066763033 & 236196 & 20476\cr
35 &  33599     &  11  &   3.3581994945 & 354294 & 28226\cr
36 &  45197     &  12  &   3.3585893055 & 531441 & 39287\cr
37 &  56039     &  14  &   3.3840009256 & 708588 & 54817\cr
38 &  81647     &  14  &   3.3598108962 & 1062882 & 75619\cr
39 &  98999     &  16  &   3.3904596729 & 1594323 & 105584\cr
40 & 163259     &  14  &   3.3324743393 & 2125764 & 146910\cr
41 &  203999    & 16  &    3.3535444722 & 3188646 & 203294\cr
42 &  241883    &  20  &   3.3881324998 & 4782969 & 283764\cr
43 &  371447    &  19  &   3.3527842988 & 6377292 & 394437\cr
44 &  540539    &  18  &   3.3332520048 & 9565938 & 547485\cr
45 &  590399    &  24  &   3.3863730003 & 14348907 & 763821\cr
46 &  907199    &  23  &   3.3532298662 & 19131876 & 1061367\cr
47 &  1081079   &  28  &   3.3828841470 & 28697814 & 1476067\cr
48 &  1851119    &  23  &  3.3261034748 & 43046721 & 2057708\cr
49 &  2041199    &  30  &  3.3725540867 & 57395628 & 2861449\cr
50 &  3243239    &  28  &  3.3350935780 & 86093442 & 3982054\cr
51 &  3840479    &  34  &  3.3638703158 & 129140163 &5552628\cr
52 &  6562079    &  28  &  3.3127733211 & 172186884 & 7721319\cr
53 &  8206559    &  33  &  3.3290528266 & 258280326 & 10758388\cr
54 &  11696759   &  33  &  3.3180085674 & 387420489 & 14994291\cr
55 &  14648759   &  38  &  3.3333603679 & 516560652 & 20866891\cr
56 &  22312799   &  36  &  3.3095614199 & 774840978 & 29079672\cr
57 &  27494879   &  42  &  3.3275907432 & 1162261467 & \cr
58 &  41746319   &  40  &  3.3053853809 & 1549681956 & \cr
59 &  \textbf{52252199}    & 46  & 3.3199050612 & 2324522934 & \cr
60 &  \textbf{78331679}    & 45  & 3.3009723129 & 3486784401 & \cr
61 &  \textbf{108606959}   & 46  & 3.2967188492 & 4649045868 & \cr
62 &  \textbf{142990559}   & 51  & 3.3016852310 & 6973568802 & \cr
63 &  \textbf{203098319}   & 52  & 3.2933942627 & 10460353203 & \cr
64 & \textbf{273985919}    & 55  & 3.2941149607 & 13947137604 & \cr
65 &  \textbf{382021919}   & 57  & 3.2893091281 & 20920706406 & \cr
66 &  \textbf{495437039}   & 63  & 3.2965467292 & 31381059609 &\cr
67 &  \textbf{681327359}   & 66  & 3.2940742853 & 41841412812 &\cr
\hline
\end{tabular}
\end{center}
}

\section{Some questions solved and proposed.}

One of the facts that our extended computation has revealed is that
sometimes the minimum in equation \eqref{defCompl} is assumed \emph{only
by the sums} and with a $j>1$.  In the range $n\le 905\,000\,000$
there are only two such instances.

The first case is the prime number $p=353\,942\,783$ (with $j=6$).
Indeed, the representation
\begin{displaymath}
353\,942\,783=2*3 + (1 + 2^2*3^2)*(2 + 3^4(1 + 2*3^{10}))
\end{displaymath}
proves that $\Vert p\Vert\le 63$, and one may verify that $\Vert
p\Vert=63$ and $\Vert p-1\Vert=63$, so that
\begin{displaymath}
\Vert p\Vert= \Vert 6\Vert+\Vert p-6\Vert=5+58=63<64=\Vert
p-1\Vert+1.
\end{displaymath}
In this case we thus have  $\Vert p\Vert=\Vert k\Vert+\Vert
p-k\Vert$ with $k=6$ (and no other choice of $k$ is adequate).

The second example is the number $n=516\,743\,639$. It is the
product of two primes $n=353\cdot1\,463\,863$. We have
\begin{displaymath}
516\,743\,639=2*3+(1+2^2 3^6)(2+3^{11})
\end{displaymath}
so that  $\Vert n\Vert\le 63$.   Also $\Vert n-1\Vert=63$, $\Vert
353\Vert =19$, $\Vert1463863\Vert=45$, $\Vert n-6\Vert=58$ and
finally $\Vert n\Vert=63$, so  that
\begin{displaymath}
1+\Vert n-1\Vert=\Vert353\Vert+\Vert1463863\Vert= 64>\Vert
6\Vert+\Vert n-6\Vert=\Vert n\Vert =63.
\end{displaymath}
Hence $\Vert n\Vert=\Vert k\Vert+\Vert n-k\Vert$ with $k=6$ and no
other choice of $k$ is adequate, as claimed.

Now we are sufficiently prepared to answer some questions asked by
Guy.

\subsection{Answering some questions of Guy}

\emph{Q1: For which values $a$ and $b$ is $\Vert 2^a3^b\Vert=2a+3b$
?}
\medskip

A1: $\Vert 2^a3^b\Vert=2a+3b$ for all $2^a3^b\le 905\,000\,000$. No
counter examples are known (to us).
\bigskip

\noindent \emph{Q2: Is it always true that  $\Vert p\Vert=1+\Vert
p-1\Vert$, if $p$ is prime ?}
\medskip

A2: No. \newline The first prime for which this is not true is $p=
353\, 942\, 783$ with $\Vert p\Vert=63$ and  $\Vert p-1\Vert=63$.
This is the only example in the range $n\le  905\, 000\, 000$.
\bigskip

\noindent \emph{Q3: Is it always true that  $3+\Vert p\Vert\le
1+\Vert 3p-1\Vert$, if $p$ is prime ?}
\medskip

A3: No.\newline There are many exceptions: $p=107$, $347$, $383$,
$467$, $587$, $683$, $719$, $887$, \dots
\bigskip

\noindent\emph{Q4: Is it always true that
 $\Vert2p\Vert=\min\bigl\{2+\Vert p\Vert,\;1+\Vert 2p-1\Vert\bigr\}$, if $p$ is prime~?}
\medskip

A4: Yes for $2p\le 905\,000\,000$. \newline \noindent Putting
$L=2+\Vert p\Vert$ and $R=1+\Vert2p-1\Vert$, we found in this range
\begin{align*}
&\Vert 2p\Vert=L\,\, ( < R )  & \text{in } 12\,317\,371 \text{ cases}\phantom{.}\\
&\Vert 2p\Vert=R\,\, ( < L ) & \text{in } \phantom{0}3\,629\,305 \text{ cases}\phantom{.}\\
&\Vert 2p\Vert=L=R  & \text{in } \phantom{0}8\,031\,758%
\text{ cases}.
\end{align*}
Note that ``$(L<R)+(R<L)+(L=R)$''$=23\,978\,434=\pi(905\,000\,000/2)$ where $\pi(\cdot)$ is the prime counting function.
\bigskip

\noindent \emph{Q5: When the value of $\Vert n\Vert$ is of the form
$\Vert a\Vert+\Vert b\Vert$, with $a+b=n$, and this minimum is not
achieved as a product, is either $a$ or $b$ equal to $1$ ?}
\medskip

A5: No.\newline We have only our two earlier mentioned ( counter )
examples: The prime $p=353942783$ and   $n=516743639$ with prime
factorization $n=353\cdot1463863$.

We have also searched in the range $n\le 905\,000\,000$ for those
cases where the minimum of $\Vert k\Vert+\Vert n-k\Vert$ is  not
assumed for $k=1$. In the cases with $k>1$ we mostly have   $k=6$,
but sometimes we have $k=8$. In all cases $\Vert n\Vert=\Vert
k\Vert+\Vert n-k\Vert=\Vert n-1\Vert$. All cases found with $k>1$
are (those with $k=8$ in boldface)
\begin{gather*}
\scriptstyle  21080618,\; 63241604,\; 67139098,\; 116385658,\;
117448688,\; 126483083,\; 152523860,\; 189724562,\;\\
\scriptstyle 212400458,\;229762259,\; 318689258,\; 348330652,\;
353942783,\; 366873514,\; 373603732,\; 379448999,\;  \\
\scriptstyle \textbf{385159320},\;
404764540,\;409108300,\;460759642,\;\textbf{465722100},\; 477258719,\; 498197068,\; 511069678,\;
\\
\scriptstyle 516743639,\; 519835084,\; 538858312,\; 545438698
,\; 545790940,\; 546853138,\; 574842670,\; \textbf{575550972},\;
 \\
\scriptstyle 581106238,\; 590785918 ,\; 608504399,\; 612752632,\;
\textbf{612752634},\; 613028608,\; 613175855,\; 614416318,\; \\
\scriptstyle 636135035,\; 637198964,\; 669796594,\; 673335934,\;
690342298,\; \textbf{690342300},\; 691406048,\; 692981240,\;\\
 \scriptstyle  \textbf{698494572},\;817595279,\; 822093928,\; 833714854,\; 860101032
,\; 861764920,\; \textbf{865717578}.
\end{gather*}
\bigskip

\noindent \emph{Q6: There are two conflicting conjectures:
\begin{displaymath}
\text{ For large $n$, }\quad (3+\varepsilon)\frac{\log n}{\log 3}
\quad\text{  ones suffice ?}
\end{displaymath}
 There are infinitely many $n$, perhaps a set of positive density for which
\begin{displaymath}
(3+c)\frac{\log n}{\log 3}\quad \text{  ones are needed, for some
$c>0$ ?}
\end{displaymath}
}
\medskip

A6: To the first question: In view of the values of  $\Vert
n\Vert/\log n$ in Table 1, the answer will most probably be no.
\medskip

A6: To the second question:  Here the answer might very well be yes.
If we solve for $c$ in the equation
\begin{displaymath}
\Vert n\Vert=(3+c)\frac{\log n}{\log 3}
\end{displaymath}
we get a mean value $\overline{c}>0.366$ and a standard deviation
$\sigma<0.047$ in the range $2\le n\le 905\,000\,000$. Also, the
frequency of the event $c>0.3$ is $>92.5\%$.

Certainly $\liminf_n\Vert n\Vert/\log n=3/\log3\approx 2.73072$. Our
computations suggest that $\limsup_n\Vert n\Vert/\log n\le 3.58$ and
that
\begin{displaymath}
\lim_{N\to\infty} \frac{1}{N-1}\sum_{k=2}^N \frac{\Vert k\Vert}{\log
k}>3 \quad\text{(possibly even $>3.06$).}
\end{displaymath}

\subsection{Some other questions}

Note that the sequence $\Vert n\Vert$ is not monotonic. It is clear
that $\Vert n-1\Vert-\Vert n\Vert\ge-1$. So, one  may pose the
question: How large can $\Vert n-1\Vert-\Vert n\Vert$ be ? We found
the first values of $n$ for which this difference is equal to $k$

\smallskip
{\scriptsize
\begin{center}
Large values of $\Vert n-1\Vert-\Vert n\Vert$
\medskip

\begin{tabular}{|r|r|r|r|r|r|r|r|r|r|}
\hline
$n$ & 6 & 12 & 24 & 108 & 720 & 1440 & 81648 & 2041200 &612360000   \\
\hline
$k=\Vert n-1\Vert-\Vert n\Vert$ & 0 & 1 & 2 & 3 & 4 & 5 & 6 & 7 & 8\\
\hline
\end{tabular}
\end{center}
} In the range  $n\le 905\, 000\, 000$  there are no  larger values
of  $\Vert n-1\Vert-\Vert n\Vert$.

\begin{conjecture}
$\displaystyle{\limsup_{n\to\infty}\;\; (\Vert n-1\Vert-\Vert
n\Vert)\;=\;+\infty.}$
\end{conjecture}

\bigskip

Let $n=\prod p_j^{a_j}$ be the standard prime-factorization of $n$.
It is clear that $\Vert n\Vert\le \sum_j \Vert p^{a_j}\Vert$. So we
define a function $\AddExc(n)= \sum_j \Vert p^{a_j}\Vert-\Vert
n\Vert$ (Additive Excess) and ask how large  $\AddExc(n)$ can be. We
found

{\scriptsize
\begin{center}
First $n_k$ with $\AddExc(n_k)=k$
\medskip

\begin{tabular}{|r|r|r|r|r|r|r|r|}
\hline $\AddExc(n)$ & 1 & 2 & 3 & 4 & 5 & 6 & 7 \\
\hline $n$ & 46 & 253 & 649 & 6049 & 69989 & 166213 &
551137 \\
\hline
$\AddExc(n)$ &  8 & 9 & 10 & 11 & & &  \\
\hline $n$ & 9064261 & 68444596 & 347562415 & 612220081 & & & \\

\hline
\end{tabular}
\end{center}
} \noindent and there are no  $n\le 905\,000\,000$ with a larger
Additive Excess.

Suppose that in our program for   $\Vert n\Vert$   we start with
$\Vert 1\Vert = 1$  and  $\Vert 2\Vert = 1+ x $  ( where $x$ is any
given real value ). Then the  $\Vert n\Vert$  will be functions of
$x$. What can be said about the resulting $\Vert n\Vert_x$ ?

\bigskip

\bigskip

Is it true that $\Vert p^k\Vert=k\Vert p\Vert$ ? Yes for $p=3$, but
we have some doubts about  $p=2$. See Section \ref{section}. We
conjecture: \emph{False for all other primes}. Examples:
\begin{align*}
\Vert 5^6\Vert&=\phantom{000\,}\Vert15\,625\Vert=29<30=6\Vert 5\Vert;\\
\Vert 7^9\Vert&=\Vert40\,353\,607\Vert=53<54=9\Vert 7\Vert;\\
\Vert 19^6\Vert&=\Vert47\,045\,881\Vert=53<54=6\Vert
19\Vert;\\
\Vert 37^5\Vert&=\Vert69\,343\,957\Vert=54<55=5\Vert 37\Vert.
\end{align*}

Our computations have revealed that for all primes $5 \le p \le 113$
($\,$with the possible exceptions $p = 73$, $97$ and $109$$\,$) it is
not true that $\Vert p^n \Vert = n \Vert p \Vert$ for all $n \ge1$.

\bigskip

We also wondered how often  $\Vert \prod  p ^ e\Vert  = \sum e \Vert
p\Vert$. We got the impression that in the long run we have about
equally often true and  false.

Pegg \cite{Pegg} asks  what  the smallest number requiring 100 ones
is? The points $(\Vert n\Vert,\log F_n)$  form approximately a
straight line (similarly as the Mahler-Popken-Selfridge points
$(m,\log M_m)$). Various least squares fits of the form $A + B t$
suggest that $M_{100}$ should be situated between
\begin{displaymath}
11\,857\,300\,000\,000\quad\text{and}\quad 27\,345\,300\,000\,000.
\end{displaymath}
A real challenge for a supercomputer$\,$! The largest number
requiring $100$ ones is $M_{100} = 7\, 412\, 080\, 755\, 407\, 364$.

Some other predictions are
\begin{align*}
F(68)&\approx 0.98 \cdot 10^9,& F(69)&\approx 1.35 \cdot 10^9,&
F(70)&\approx 1.86 \cdot 10^9,\\
F(71)&\approx 2.56 \cdot 10^9,& F(72)&\approx 3.53 \cdot 10^9,&
F(73)&\approx 4.85 \cdot 10^9,\\
F(74)&\approx 6.68 \cdot 10^9,& F(75)&\approx 9.20 \cdot 10^9,&
F(80)&\approx 45.54 \cdot 10^9.
\end{align*}

\subsection{Is it true that $\Vert 2^k\Vert=2k$~?}\label{section}

Selfridge asked whether $\Vert 2^k\Vert=2k$ for all $k\ge1$. We have
verified this for $1\le k\le29$. Nevertheless, we will present an
argument suggesting  that the answer may very well be no.

Given a natural number $n$ with complexity $\Vert n\Vert=a$ we
denote by $M_a$ the greatest number with the same complexity, and we
will call
\begin{displaymath}
\CR(n)=1-\frac{n}{M_a}
\end{displaymath}
the \emph{complexity ratio} of $n$.

We  always have $0\le CR(n)< 1$. In a certain sense the numbers $n$
with a small complexity ratio are \emph{simple} and those with a
large complexity ratio are \emph{complex}.  To illustrate this we
present here some numbers comparable in size but with different
complexity ratios and their corresponding minimal representations.

\begin{center}
\begin{tabular}{|r|r|l|l|}
\hline
$n$ &  $\scriptstyle{\Vert n\Vert}$ & $\scriptstyle{\CR(n)}$ &  \text{Minimal Expression} \\
\hline
371447 & 43 & 0.94 & $\scriptstyle{1+2(1+2(1+2^2(1+2^2)(1+2(1+2^4(1+2^43^2)))))}$ \\
373714 & 40 & 0.82 & $\scriptstyle{2(1+2^3(1+2^2(1+2\cdot3(1+2^23^5))))}$ \\
377202 & 39 & 0.76 & $\scriptstyle{3(1+2\cdot3)^2(1+(1+2^2)(1+2\cdot3^2)3^3)}$ \\
377233 & 38 & 0.65 & $\scriptstyle{(1+2^5 3)(1+2^4 3^5)}$ \\
360910 & 37 & 0.49 & $\scriptstyle{(1+2\cdot 3^3)(1+3^8)}$ \\
422820 & 37 & 0.40 & $\scriptstyle{2^2(1+2^4 3^2)3^6}$ \\
492075 & 37 & 0.31 & $\scriptstyle{(1+2^2)^2 3^9}$ \\
413343 & 36 & 0.22 & $ \scriptstyle{(1+2\cdot3)3^{10}} $ \\
531441 & 36 & 0   & $\scriptstyle{3^{12}} $ \\
\hline
\end{tabular}
\end{center}

Let $S$ be some (arbitrary but fixed) natural number (this will be
the span of $n$).  Choose $S$ not too small. For example,
$S=1\,000\,000$.
\medskip
Let $M=\max\{\Vert s\Vert:1\le s\le S\}$. So, $M$ is also fixed.
\medskip

Now choose $k$ such that $2^{3k}>S$. Clearly there are infinitely
many such $k$.

Now let $n$ satisfy $2^{3k}-S\le n<2^{3k}$, and let $\Vert
n\Vert=3a+r$ with $0\le r\le2$.

Then we have
\begin{displaymath}
\CR(n)=1-\frac{n}{3^a},\quad 1-\frac{n}{4\cdot 3^{a-1}},\quad
1-\frac{n}{2\cdot 3^{a}}
\end{displaymath}
for $r=0$, $1$ or $2$, respectively. Therefore, in all cases we will
have
\begin{displaymath}
\CR(n)=1-f\frac{n}{3^a}
\end{displaymath}
where $f=1$ for $r=0$, $f=3/4$ for $r=1$ and $f=1/2$ for $r=2$.
\medskip

Now choose a small fixed $p>0$, ($p=1/1000$, say).

Let's now consider the inequality
\begin{equation}\label{E}
\CR(n)+p<1-\Bigl(\frac89\Bigr)^k.
\end{equation}
For large $k$ this comes very close to the event $\CR(n)+p\le1$ or
$\CR(n)\le1-p$. Quite extensive statistics on $\CR(n)$ suggest
strongly that this event is highly probable (for small $p>0$). So,
we venture to \emph{assume} that we have \eqref{E}. Observe that
this is equivalent to
\begin{displaymath}
\Bigl(1-f\frac{n}{3^a}\Bigr)+p\le 1-\Bigl(\frac89\Bigr)^k.
\end{displaymath}
Hence, since $f\le1$ (also using previous assumptions)
\begin{displaymath}
\frac{2^{3k}}{3^a}>\frac{n}{3^a}\ge f \frac{n}{3^a} \ge
p+\Bigl(\frac89\Bigr)^k=p+\frac{2^{3k}}{3^{2k}}
\end{displaymath}
so that $2k>a$ or $2k-a>0$.

Also observe that
\begin{displaymath}
p+\Bigl(\frac89\Bigr)^k\le 1-\CR(n)=f\frac{n}{3^a} \le
3^{2k-a}\frac{n}{3^{2k}}<3^{2k-a}\frac{2^{3k}}{3^{2k}}=3^{2k-a}\Bigl(\frac89\Bigr)^k
\end{displaymath}
so that
\begin{displaymath}
p+\Bigl(\frac89\Bigr)^k<
3^{2k-a}\Bigl(\frac89\Bigr)^k\quad\text{or}\quad
p\Bigr(\frac98\Bigr)^k+1<3^{2k-a}.
\end{displaymath}
Now choose $k$ so large that $3^{M+r}\le
3^{M+2}<p\left(\frac98\right)^k+1$, without violating previous
assumptions.

Then we  clearly have $2k-a>M+r$.

Now we can conclude that
\begin{multline*}
\Vert 2^{3k}\Vert=\Vert n+(2^{3k}-n)\Vert\le\Vert n\Vert+
\Vert 2^{3k}-n\Vert\le 3a+r+\Vert\text{some } s\le S\Vert\le\\
\le 3a+r+M<3a+(2k-a)=2a+2k<2(2k)+2k=6k
\end{multline*}
so that
\begin{displaymath}
\Vert 2^{3k}\Vert<6k=3k\Vert 2\Vert.
\end{displaymath}

Hence, the answer to Selfridge's question might very well be no.

\begin{center}
\begin{figure}[H]
  \includegraphics[width=10cm]{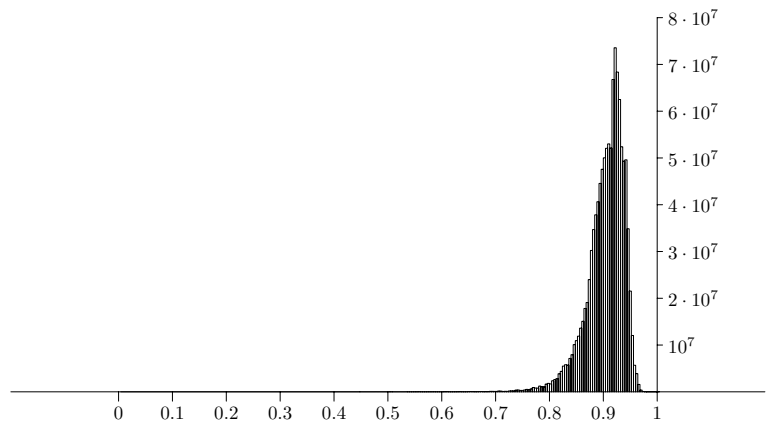}\\
  \caption{Distribution of $\CR(n)$ for $1\le n\le 905\,000\,000$}\label{figura1}
\end{figure}
\end{center}

\vspace{1cm}

\noindent \textsc{J. Arias de Reyna, Facultad de Matemáticas,
Universidad de Sevilla,
 Apdo.\-~1160, 41080-Sevilla, Spain.}
e-mail {\tt arias@us.es}.\hfil\break Supported by  grant
MTM2006-05622.
\bigskip

\noindent \textsc{J. van de Lune, Langebuorren 49, 9074 CH Hallum,
The Netherlands} (formerly at CWI, Amsterdam).  e-mail {\tt
j.vandelune@hccnet.nl}

\end{document}